\documentclass[a4paper,11pt]{amsart}

\usepackage{amsthm}
\usepackage{amsmath}
\usepackage{amsfonts}
\usepackage{amssymb}
\usepackage{times}

\theoremstyle{definition}
\newtheorem{defn}{\indent\bf Definition}
\newtheorem{rem}[defn]{\indent\bf Remark}

\theoremstyle{plain}

\newtheorem{thm}[defn]{\indent\bf Theorem}

\begin{document}

\title[Hecke algebra of double cosets]{A universal, non-commutative $C^*$-algebra\\ associated to the Hecke algebra of double cosets}
\author[Florin R\u adulescu]{Florin R\u adulescu${}^*$
\\ \\ 
Dipartimento di Matematica\\ Universita degli Studi di Roma ``Tor Vergata''}

\maketitle 

\thispagestyle{empty}
\renewcommand{\thefootnote}{}
\footnotetext{${}^*$ Member of the Institute of  Mathematics ``S. Stoilow" of the Romanian Academy}

\def\tilde{\widetilde}
\def\a{\alpha}
\def\T{\theta}
\def\PSL{\mathop{\rm PSL}\nolimits}
\def\SL{\mathop{\rm SL}\nolimits}
\def\PGL{\mathop{\rm PGL}}
\def\Per{\mathop{\rm Per}}
\def\GL{\mathop{\rm GL}}
\def\Out{\mathop{\rm Out}}
\def\Int{\mathop{\rm Int}}
\def\Aut{\mathop{\rm Aut}}
\def\ind{\mathop{\rm ind}}
\def\card{\mathop{\rm card}}
\def\d{{\rm d}}
\def\Z{\mathbb Z}
\def\R{\mathbb R}
\def\Q{\mathbb Q}
\def\N{\mathbb N}
\def\C{\mathbb C}
\def\bH{\mathbb H}
\def\L{{\mathcal L}}
\def\G{{\mathcal G}}
\def\U{{\mathcal U}}
\def\H{{\mathcal H}}
\def\A{{\mathcal A}}
\def\S{{\mathcal S}}
\def\O{{\mathcal O}}
\def\D{{\mathcal D}}
\def\B{{\mathcal B}}
\def\K{{\mathcal K}}
\def\cC{{\mathcal C}}
\def\ptimes{\mathop{\boxtimes}\limits}
\def\potimes{\mathop{\otimes}\limits}

\begin{abstract}
Let G be a discrete group and $\Gamma$ an almost normal subgroup. The operation of cosets concatanation extended by linearity gives rise to an operator system that is embeddable in a natural C* algebra. The Hecke algebra naturally embeds as a diagonal of the tensor product of this algebra with its opposite. When represented on the $l^2$ space of the group, by left and right convolution operators, this representation gives rise to abstract Hecke operators that in the modular group case, are unitarily equivalent to the classical operators on Maass wave forms\end{abstract}

\section*{}
Let $G$ be a discrete group, $\Gamma$ an almost normal subgroup.
In this work, we introduce a formalism for the free algebra having as free generators the cosets of $\Gamma$ in $G$, subject to the
relations that define the Hecke algebra and its standard action on the space of cosets.

Let  $\C(\Gamma\backslash G)$, respectively  $\C(G\backslash\Gamma)$ be the $\C$-linear space having as a linear basis the  left and respectively right cosets. Let 
$\C(\Gamma\backslash G/G)$ be the $\C$-linear space  having as a linear basis the  double cosets of $\Gamma$ in $G$. Since $\Gamma$ is almost normal, it follows that
every double coset $[\Gamma\sigma\Gamma]$ is a finite union (of which we think as of a sum in the linear space
$\C(\Gamma\backslash G)$, respectively $\C(G\backslash\Gamma)$ ) of left  and respectively right cosets of $\Gamma$ in $G$.

For $\sigma$ in $G$, let $\Gamma_\sigma = \sigma\Gamma \sigma^{-1} \cap \Gamma$. This subgroup is, by hypothesis, 
of finite index in $G$. Let  $\S$ be the lattice  of  subgroups  generated by the above subgroups. Note that all the subgroups in $\S$ have finite index.

We endow the Hilbert space $\ell^2(\Gamma\backslash G)$ with the canonical Hilbert space structure  such that the cosets
$[\Gamma\sigma]$, $\sigma \in G$, where $\sigma \in G$ runs over a system representatives of right cosets of 
$\Gamma$ in $G$, is an orthonormal basis ((and similarly for
$\ell^2(G \backslash \Gamma)$ .

There is a  canonical      prehilbert space structure on $\bigcup_{\Gamma_0\in \S} \ell^2(\Gamma_0/G)$ that  is uniquely determined by the requirement 
 that $[\sigma_1\Gamma] = \sum[\sigma_1 s_i\Gamma_0]$, for
all $\sigma_1\in G$ where $s_i$ are a system of representatives for right cosets
of $\Gamma_0$ in $\S$ and by the  requirement   that all the inclusions $\ell^2(\Gamma/G)$ in $\ell^2(\Gamma_0/G)$, for
all $ \Gamma_0\in \S$
are isometric.

The space  $\C(\Gamma\backslash G/\Gamma)$ has a natural action on left and respectively right cosets. The action of double coset $[\Gamma\sigma_1\Gamma]$ on $[\Gamma\sigma_2]$ is simply
the projection of $[\sigma_1\Gamma \sigma]$, which, by linearity,  belongs to the reunion $\bigcup_{\Gamma_{\sigma_2}\in \S} \ell^2(\Gamma_{\sigma_2})$
onto $\ell^2(\Gamma/G)$. Naturally, we have a similar left action of $\C(\Gamma\backslash G/\Gamma)$
on $\ell^2(G/\Gamma)$.  This actions define canonically a $\C-$ algebra structure on $\C(\Gamma\backslash G/\Gamma)$, and with this structure the space $\C(\Gamma\backslash G/\Gamma)$ is the Hecke algebra associated with the inclusion $\Gamma$ in $G$.

The relations defining the Hecke algebra are thus summarized by the requirement that 
$$
\sum_i[\sigma_1^i\Gamma][\Gamma\sigma_2^i]=
\sum_j[\theta_1^j\Gamma][\Gamma\theta_2^j]\eqno(*)
$$
if the unions $\bigcup_i  \sigma_1^i\Gamma \sigma_2^i$, $\bigcup_j \theta_1^j\Gamma \theta_2^j$ are
disjoint and equal. In particular $[\Gamma]$ is the unit.

This allows, naturally, to define an abstract $*$-algebra, whose  free generators are the cosets
(left or right) of $\Gamma$ in $G$, subject to the relations $(*)$. By analogy with  the Jones's algebra of  higher relative
commutants ([Jo], [Bi]) (or the algebra  of bimodules over $\Gamma$ in $\ell^2(G)$), we define a larger algebra
$\B(\Gamma,G)$ which in addition allows to have a consistent definition for the  localization of the support of the cosets (the support is the coset itself viewed as a characteristic function).
This is done by adding to the algebra of cosets the algebra of characteristic functions of cosets of groups
in $\S$.

Our main result is that, under the assumption that there exists a (projective) unitary representation of $G$ on
$\ell^2(\Gamma)$ that extends the left regular representation (with the cocycle induced from $\pi$) on
$\ell^2(\Gamma)$, the algebra $\B(\Gamma,G)$ admits a unital $C^*$-representation, which is faithful on $\C(\Gamma\backslash G/\Gamma )$ and is consistent with its action on 
$\C(G\backslash\Gamma)$. (Examples of such representations
are found in [GHJ]).
Here the * operation is obtained as a linear extension of $[\Gamma\sigma]^*=[\sigma^{-1}\Gamma]$, $\sigma$ in $G$.

The construction of the $C^*$-representation of the algebra $\B(\Gamma,G)$ is based on the analysis ([Ra]) of
the properties of the positive definite function on $G$
$$
\varphi_I (g) = \langle \pi(g)I,I\rangle,  g\in G
$$
where $I \in \ell^2(\Gamma)$ is the vector corresponding to the identity of $\Gamma$ in $\ell^2(\Gamma),$.

We let $X_\Gamma$ be the (commutative) subalgebra of $\ell^{\infty} (G)$ generated by characteristic functions of  cosets of groups in $\S$.

Let $\chi_{\Gamma} \in \ell^{\infty} (G)$   be the characteristic function of the group $\Gamma$. The unit of the algebra $\B(\Gamma,G)$ will be  $\chi_{\Gamma}$.  We will  prove that there exists a canonical $*$-representation of the Hecke algebra 
$\C(\Gamma\backslash G/\Gamma)$ into $\B(\Gamma,G)\bigotimes\limits_{X_{\Gamma}} \B(\Gamma,G)^{\rm op}$,
defined by the mapping
$$\C(\Gamma\backslash G/\Gamma) \ni [\Gamma\sigma\Gamma]\mapsto \chi_\Gamma([\Gamma\sigma\Gamma]\otimes
[\Gamma\sigma\Gamma]^*) \chi_\Gamma, \sigma \in G.$$

This map is, in  the canonical representation of 
$\B(\Gamma,G)\bigotimes\limits_{X_{\Gamma}} \B(\Gamma,G)$ as 
left and right convolution operators on $\ell^2(\Gamma)$,  unitarily equivalent to classical
Hecke operators on Maass forms ([Ra]).

We extend this representation to a representation of $\C(\Gamma\backslash G)$ into
$\B(G/\Gamma)\bigotimes\limits_{X_{\Gamma}}$  $\B(\Gamma/G)^{\rm op}$
that is compatible with the action of the Hecke algebra $\C(\Gamma\backslash G/\Gamma)$ onto
$\C(\Gamma/ G)$. 

We remark that this construction gives a canonical operator system structure on the linear space 
$\{\C([\sigma_1\Gamma\sigma_2])\mid \sigma_1, \sigma_2 \in G\}$ having as basis the sets $[\sigma_1\Gamma\sigma_2])\mid \sigma_1, \sigma_2 \in G$.  This space is  is isomorphic to
$
\C(G\backslash\Gamma) \bigotimes\limits_{\C(\Gamma\backslash G/G)} \C(\Gamma/ G)$.

This last space is canonically included in the similar spaces, constructed by using subgroups of the form
$\Gamma_\sigma$ instead of $\Gamma$. It is conceivable that one could obtain such an operation system
structure, canonically associated to the Hecke algebra, simultaneously for all levels ($\Gamma_\sigma$ in $\S$).

We start with the definition of the free algebra of cosets of $G$ in $\Gamma$.

\begin{defn}
Let $\Gamma$ be an almost normal subgroup of a discrete group $G$.

Let $\C(\Gamma\backslash G)$, $\C(G\backslash\Gamma)$, $\C(\Gamma\backslash G/\Gamma)$ be
the linear vector space having as basis, respectively the left, right, and
double cosets of $\Gamma$ in $G$.
For $\sigma$ in $G$ we denote  such a coset (respectively a double coset)  by   
$[\Gamma\sigma]$, $[\sigma\Gamma]$ and by $[\Gamma\sigma\Gamma]$ respectively.

Let $I(G,\Gamma)$ be the free $\C$-algebra whose (algebra) generators are  
$[\Gamma\sigma]$, $[\sigma\Gamma]$,  for $\sigma$ in $G$. We define a natural $*$-operation
by requiring $[\Gamma\sigma]^*=[\sigma^{-1}\Gamma]$.

Consider $J(G,\Gamma)$ the double sided ideal corresponding to all relations of the form 
$$
\sum_i[\sigma_1^i\Gamma][\Gamma\sigma_2^i]=
\sum_j[\theta_1^j\Gamma][\Gamma\theta_2^j]\eqno(1)
$$
whenever $\sigma_{\alpha}^i$, $\theta_{\beta}^j$, $\alpha=1,2$, $\beta=1,2$, $i=1,2,\ldots, N$, $j=1,2,\ldots,M$
are a finite set of elements in $G$ such that the unions
$\bigcup_i [\sigma_1^i\Gamma\sigma_2^i]$,  $\bigcup_j [\theta_1^j\Gamma\theta_2^j]$
are disjoint, and equal. We also assume the relation that $[\Gamma]$ is the identity of the algebra.

Then the universal $*$-algebra associated to $\Gamma$, $G$ will be $$\A(\Gamma, G)=I(\Gamma,G)/J(\Gamma,G).$$
Then the Hecke algebra  $\C(\Gamma\backslash G/\Gamma)$ is a subalgebra of  
$\A(\Gamma, G)$, and the embedding is compatible with the left (and right) action of 
$\C(\Gamma\backslash G/\Gamma)$ on $\C(\Gamma/ G)$ (and respectively $\C(G\backslash \Gamma)$).
\end{defn}

{\it Proof.}
Indeed the relations defining the action of the Hecke algebra, on the space of left and right
cosets, and its multiplication are of the form of the relation in~(1).

For example, the action of a double coset $[\Gamma\sigma_1\Gamma]$ on a coset 
$[\Gamma\sigma_2]$ is 
$$
[\Gamma\sigma_1\Gamma][\Gamma\sigma_2]=\sum_j[\Gamma][\Gamma\sigma_1r_j\sigma_2]
$$
if $r_j$ are representatives for $\Gamma_{\sigma_1}=\sigma_1\Gamma\sigma_1^{-1}\cap\Gamma$
(that is $[\Gamma\sigma_1\Gamma]=\bigcup[\Gamma\sigma_1r_j]$).

On the other hand, if $s_j$ are left representatives that is, 
$[\Gamma\sigma\Gamma]=\bigcup[s_i\sigma\Gamma]$)
then we get
$$
\sum[s_1\sigma\Gamma][\Gamma\sigma_2]=\sum_j[\Gamma][\Gamma\sigma_1r_j\sigma_2].\quad\Box
$$
\vskip4pt

We will prove in the sequel that the algebra $\A(\Gamma, G)$ admits a $*$-representation into a
$C^*$-algebra, and thus that there exists a maximal $C^*$-algebra associated with $\A(\Gamma, G)$.

\begin{rem}
Let $\theta$ be an automorphism of $G$ preserving $\Gamma$. Then $\theta$ obviously extends
to an automorphism of $\A(\Gamma, G)$.

\

\

We also have a simple method to algebraically describe the localization of the supports of 
$[\Gamma\sigma]$, $[\sigma\Gamma]$, $\sigma \in G$, by extending $C^*$-algebra. For this purpose we  consider a larger $*$-algebra.

Let $X_\Gamma$ be the $\C$-subalgebra of $L^\infty(G)$ consisting of characteristic functions of left and
right cosets in $G$ of subgroups in $\S$ (with $/S$ as above).
\end{rem}

\begin{defn}
The localized $*$-free $\C$-algebra of cosets of $\Gamma$ in $G$. We consider this time the free $\C$-algebra
$I_1(\Gamma,G)$ whose generators are cosets $[s_1 \Gamma_ {\sigma_1}]$,  $[\Gamma_{ \sigma_2} s_2]$,
where $\sigma_1, \sigma_2 \in G$, $s_1, s_2$ in $G$ and the algebra $X_\Gamma$. Let $J_1(\Gamma,G)$
be the bilateral ideal in $I_1(\Gamma,G)$ corresponding to the following relations.

(0) If $C$ is a coset of some modular subgroup in $\S$, which in turn is a disjoint reunion $\bigcup_j D_j$
of cosets then
$$
[C] = \sum [D_j].
$$

(1) Fix two cosets $C_1$, $C_2$ (left or right) of subgroups in $\S$
(that is, $C_i=[\sigma_i\Gamma_{\theta_i}]$ or $\Gamma_{\theta_i}\sigma_i$ for some
$\sigma_i$, $\theta_i$ in $G$, $i=1,2$).

By passing to smaller cosets there exists partitions (into cosets)
$C_1=\bigcup_{i\in I} A_i$, $C_1C_2=\bigcup_{j\in J}D_j$ and a map $\pi:I\to J$
such that if $c\in A_i$ then $cC_2\in D_{\pi(i)}$.

Then we require that
$$
[C_1]\chi_{[C_2]}=\sum_i\chi_{[D_{\pi(i)}]}[A_i].
$$

For example $[\sigma\Gamma]\chi_{\Gamma}=\chi_{[\sigma\Gamma]}[\sigma\Gamma]$ and
$[\Gamma\sigma]\chi_{\Gamma}=\chi_{[\sigma\Gamma]}[\Gamma\sigma]$.
This corresponds to the fact that, when the group algebra of $G$ acts on $L^2(G)$, then
for every subset $A$ of $G$, we have $g\chi_A g^{-1}=\chi_{gA}$, for all $g$ in $G$.

This property will correspond to the fact that the quotient algebra over $J_1(\Gamma,G)$
will be, as a linear space, the linear span $Sp (A(\Gamma,G)X_\Gamma)=Sp (X_\Gamma \A(\Gamma,G))$.

The relation (1) from the previous definition becomes:

(2) for all $\sigma_1^i, \sigma_2^i$ in $G$, if $\bigcup [\sigma_1^i\Gamma \sigma_2^i]$
is a disjoint union and its further equal to the disjoint union $\bigcup C_i$ of
cosets then
$$
\sum[\sigma_1^i\Gamma][ \Gamma \sigma_2^i]=\sum_j [C_j].
$$

In addition to this property we will add some properties which hold true in our model
and are due to the fact that our algebra is related to the Jones's construction ([Jo])
for the inclusion $\Gamma\subset G$, and its Pimsner-Popa basis ([PP])

(3) Let $C_1 = [\sigma_1\Gamma]$, $C_2=[\sigma_2\Gamma]$ be two cosets of $\Gamma$ in $G$.
Then $\chi_\Gamma [\sigma_1\Gamma][ \Gamma \sigma_2]\chi_\Gamma = \delta_{C_1,C_2}\chi_\Gamma$,
where $\delta$ is the Kronecker symbol.

(4) For all $\sigma$ in $G$, if $[\Gamma \sigma\Gamma]$ is the disjoint reunion of $[s_1 \sigma\Gamma]$
then 
$$
\sum[\Gamma\sigma s_i]\chi_\Gamma[s_i\sigma\Gamma]= \chi_{[\Gamma \sigma\Gamma]}.
$$

(5) In the same conditions as in (4) we have that 
$$
\sum[\Gamma\sigma s_i][s_i\sigma\Gamma]= [\Gamma: \Gamma_\sigma]{\rm Id}.
$$

We define then the localized free algebra of cosets as the algebra
$\B(\Gamma,G) = I_1(\Gamma,G)/J_1(\Gamma,G)$. It is an algebra with generators cosets
in $G$ of subgroups in $\S$.

Our main result is that the $*$-algebra $\B(\Gamma,G)$ admits a $*$ unital representation into
a $C^*$-algebra, (more precisely into a $C^*$-subalgebra of  $\L(G)$, the II$_1$ factor associated to the discrete
group $G$).
\end{defn}

\begin{thm}
Let $\L(G)$ be the finite von Neumann algebra associated to $G$. Assume that
there exists a unitary (eventually a projective unitary representation with 2 cocycle $\varepsilon$) $\pi$
of $G$ on $\ell^2(\Gamma)$ that extends the left regular unitary representation of $\Gamma$ on 
$\ell^2(\Gamma)$ (eventually with cocycle $\varepsilon$). Then there exists an embedding of 
$\B(\Gamma,G)$ into $\L(G)$, such that every element
$[\sigma_1\Gamma][\Gamma\sigma_2]$ = $[\sigma_1\Gamma\sigma_2]$ in $\B(\Gamma,G)$ , $\sigma_1,\sigma_2\in G$
is supported in $\L(G) \cap \ell^2(\sigma_1\Gamma\sigma_2)$.

Examples of such representations can be found in {\rm ([GHJ], [Jo])}.
\end{thm}

{\it Proof.} This is based on the results proven in [Ra]. Let $u$ be
a unitary in $\L(\Gamma)$ (viewed as trace vector) and let $\varphi_J^u$ be the positive definite functional on $G$, defined
as the matrix coefficient
$$
\varphi_J^u(g) = \overline{\langle \pi(g)u,u\rangle}_{\ell^2(\Gamma)}, \quad g\in G.
$$

For simplicity when $u$ is the identity 1 we denote
$$
t(g) = \varphi_J^1 (g),  \quad g\in G.
$$

We proved in [Ra] that the matrix coefficients $t$ have the property that
$$
t(\theta_1\theta_2)=\varepsilon(\theta_1,\theta_2)
\sum_{\gamma\in\Gamma}t(\theta_1\gamma^{-1})t(\gamma\theta_2),\quad \theta_1,\theta_2\in G.
$$

We denote for $A$ a subset of $G$,
$$
t^A=\sum_{\theta\in A}t(\theta) \theta\in\L(G).
$$
Then ([Ra])
$$
t^{\sigma_1\Gamma}t^{\Gamma\sigma_2}=t^{\sigma_1\Gamma\sigma_2},\quad \sigma_1,\sigma_2\in\Gamma.
$$

We also proved in [Ra], that by a change coordinates we may assume that all $t^{\Gamma\sigma}$
belong to $\L(G)$.

To define the representation of $\B(\Gamma,G)$, we associate to each coset $C$ the element
$t^C$, while to $\chi_C$ we associate the same characteristic function of $C$, viewed as a bounded
linear operator on $\ell^2(\Gamma)$, acting by multiplication.

Note that $t^{\Gamma\sigma}$ has in fact the formula $(\pi(\sigma)1)^*\sigma$, $\sigma\in G$.

The properties (0), (1), (2) are obvious.
It remains to show the properties (3), (4), (5) in terms of the Jones's construction 
for $\L(\Gamma) \subseteq L(G)$ ([Jo]).

Property (3) is equivalent to $E_{\L(\Gamma)}^{\L(G)}((t^{\Gamma\sigma_1})^*(t^{\Gamma\sigma_2}))
= \delta_{\Gamma \sigma_1, \Gamma \sigma_2}{\rm Id}$, $\sigma_1, \sigma_2 \in G$.
But $(t^{\Gamma\sigma_1})^*(t^{\Gamma\sigma_2})= t^{\sigma_1^{-1}\Gamma\sigma_2}$, $\sigma_1, \sigma_2 \in G$
and $t^{\sigma_1^{-1}\Gamma\sigma_2}$ has a nonzero coefficient in $\Gamma$ if and only if $\sigma_1=\sigma_2$
and in this case the coefficient is exactly the identity.

It then follows that $(t^{\Gamma\sigma})\chi_\Gamma$ is a partial isometry, and if $
[\Gamma\sigma\Gamma] = \bigcup [\Gamma\sigma s_i]$, where $s_i$ are representation
then $(t^{\Gamma\sigma s_i})\chi_\Gamma$ are orthogonal partial isometries with initial space
$\ell^2(\Gamma)$ and range contained in $\ell^2(\Gamma\sigma\Gamma)$.
Because there are $[\Gamma : \Gamma_\sigma]$ such isometries, and the $\L(\Gamma)$ bimodule
$\ell^2(\Gamma\sigma\Gamma)$ has multiplicity $[\Gamma : \Gamma_\sigma]$ it follows that the sum
of the orthogonal ranges of the partial isometries $(t^{\Gamma\sigma s_i})\chi_\Gamma$ is exactly
$\ell^2(\Gamma\sigma\Gamma)$.

This proves property (4) and property (5) is as in [PP]. Note that property (5) also proves
the boundedness of $t^{\Gamma_\sigma}$. $\quad\Box$

We will now prove that the algebra $\A(\Gamma,G)$ has a special diagonal representation into 
$\B(\Gamma,G)\bigotimes\limits_{X_\Gamma} \B(\Gamma,G)^{\rm op}$. Here the tensor
product is with amalgamation over $\chi_\Gamma$ (with $\B(\Gamma,G)^{\rm op}$ acting from the
right).

For the purpose of calculations we need to consider a larger algebra containing $\B(\Gamma,G)$
by taking into account the action of $\Gamma$.

\begin{defn}
The localized $*$ $\Gamma$-free algebra $\cC(\Gamma,G)$ of cosets of modular subgroups of $\Gamma$ in $G$,
is the $*$-algebra generated by $\B(\Gamma,G)$ and the group $*$-algebra of $\Gamma$, where
the action of $\Gamma$ on $X_\Gamma$ is compatible with the action of cosets on $X_\Gamma$, that is
$\gamma\chi_A\gamma^{-1} = \chi_{\gamma A}$, $\gamma\in \Gamma$, $A$ subset of $G$.
\end{defn}

In fact, our representation of the algebra $\B(\Gamma,G)$ proves that the algebra $\L(G)$ plays
the role as an algebra of coefficients. Thus we have the following additional property:

(6) $[\Gamma \sigma]\sigma^{-1}$ belong all to $\L(\Gamma)$, $\sigma \in G$. 

We denote for a coset $C$ of a modular subgroup in $\S$, by $C_\gamma$ , 
$\gamma\in \Gamma$, the sum
$$C_\gamma=\sum_{i}\gamma D_i\gamma^{-1},$$
if  as a disjoint reunion we have that
$$C=\bigcup_i  \gamma D_i\gamma^{-1},$$
where $D_i$ is a finite family of cosets of groups in $\S$.
 As a consequence of the fact that ${\rm Ad}\, \gamma$ is an inner
automorphism of $\cC(\Gamma,G)$  we have with the notations of property (2) in the preceding definition
that 
$$
\sum_i [\sigma_1^i \Gamma]_\gamma [\Gamma\sigma_2^i]_\gamma = \sum_j [C_j]_\gamma.
$$

We now prove a diagonal representation of the Hecke algebra, which we extend to the action of the 
Hecke algebra $\C(\Gamma\backslash G/\Gamma)$ on $\C(\Gamma\backslash G)$.

Consider the algebra $\B(\Gamma,G)^{\rm op}$ where multiplication of $b_1 \odot b_2$, $b_1,b_2 \in \B(\Gamma,G)^{\rm op}$ is by definition $b_2b_1$. With this we have the following

\begin{thm}
Define the following map
$$
\Phi : \C(\Gamma\backslash G) \to \chi_\Gamma ( \B(\Gamma,G) 
\textstyle\bigotimes\limits_{X_\Gamma} \B(\Gamma,G)^{\rm op}) \chi_\Gamma
$$
by
$$
\Phi([\Gamma\sigma])=\chi_{\Gamma}[\Gamma\sigma]\otimes[\Gamma\sigma\Gamma]^*\chi_{\Gamma},\quad \sigma\in G.
$$
Then $\Phi$ is a $\C(\Gamma\backslash G/\Gamma)$-linear map, that is for all $\sigma_1,\sigma_2$ in $G$, 
$$
\Phi(\Gamma\sigma_1\Gamma) \Phi(\Gamma\sigma_2) = \Phi([\Gamma\sigma_1\Gamma][\Gamma\sigma_2]).
$$

Note that in particular this proves that the map
$$
[\Gamma\sigma\Gamma] \to \chi_\Gamma ([\Gamma\sigma\Gamma]\otimes [\Gamma\sigma\Gamma]^*)\chi_\Gamma
$$
is a $*$-algebra morphism from the Hecke algebra into 
$$
\chi_\Gamma  ( \B(\Gamma,G) 
\textstyle\bigotimes\limits_{X_\Gamma} \B(\Gamma,G)^{\rm op}) \chi_\Gamma.
$$

As proved in {\rm [Ra]}, these maps, when representing 
$\B(\Gamma,G) \textstyle\bigotimes\limits_{X_\Gamma} \B(\Gamma,G)^{\rm op}$ on $\ell^2(\Gamma)$
are unitarily equivalent to classical Hecke operators.
\end{thm}

{\it Proof.} Denote by $\B = \B(\Gamma,G)$. To make the calculations, we write
elements $X$ in $\B\textstyle\bigotimes\limits_{X_\Gamma} \B^{\rm op}$
in the form $X =\sum [C_i] \otimes [D_i] f_i$, where $[C_i]$, $[D_i]$ are cosets and $f_i$ functions
in $X_\Gamma$. For $\gamma$ in $\Gamma$ we denote by 
$$
X(\gamma) = \sum_i  [C_i] \otimes [D_i] f_i(\gamma).
$$
Working in the algebra $\cC = \cC(\Gamma,G)$ instead
of $\B=\B(\Gamma,G)$, we have that for $\gamma$ in $\Gamma$, $A,B$ cosets in $\cC(\Gamma,G)$ , and $e$ the identity of the group $\Gamma$,  
$$
(A\otimes B)f(\gamma) = (A\otimes \gamma B) f(e).
$$
This holds true since $B$ acts on the right on $X_\Gamma$. 

This formula has a meaning also if we have another function $h$ in front of $A\otimes B$. By
$h(A\otimes B)(\gamma)$ we will mean the result evaluated at $\gamma$, 
when $h$ is moved to the left.

Then $(A\otimes \gamma B) f(e)$ is further equal to
$$
(A \otimes \gamma \odot (\gamma B \gamma^{-1}))f(e) = (A \otimes (\gamma \odot B_\gamma))f(e).
$$

Note the obvious identity for all cosets $A,B;M,N$
$$
(M \otimes N)(A \otimes (\gamma \odot B_\gamma)) = 
MA \otimes(\gamma \odot N_\gamma \odot B_\gamma)=
(1\otimes \gamma) (MA \otimes( N_\gamma \odot B_\gamma)). \eqno(2)
$$
Then, if $\sigma\in G$, $[\Gamma\sigma\Gamma] = \sum [\Gamma\sigma s_i]$ with $s_i$ representatives, then
\begin{align*}
\Phi([\Gamma\sigma])(\gamma) & =\chi_{[\Gamma]}([\Gamma\sigma]\otimes([\Gamma\sigma\Gamma])^*)(\gamma)
= \chi_{[\Gamma]}  \sum_i [\Gamma\sigma]\otimes ([\Gamma\sigma s_i])^*(\gamma)
\\
& = \chi_{[\Gamma]}  \sum_i [\Gamma\sigma]\otimes (\gamma\odot[\Gamma\sigma s_i]_\gamma)(e)=
\\
=\chi_{[\Gamma]} (1\otimes \gamma)\sum_i [\Gamma\sigma]\otimes ([\Gamma\sigma s_i]_\gamma)(e)=
\\
 (1\otimes \gamma)\chi_{[\Gamma]}\sum_i [\Gamma\sigma]\otimes ([\Gamma\sigma s_i]_\gamma)(e) .
\end{align*}
Passing $\chi_{[\Gamma]}$ to the right hand side, 
this is further equal to
$$
\sum_i [\Gamma\sigma]\otimes (\gamma\odot[\Gamma\sigma s_i]^*_\gamma)\chi_{\sigma^{-1}\Gamma(\sigma s_i)^{-1}}(e).
$$

The only non zero terms that remains in the sum is 
$$\Phi([\Gamma\sigma])(\gamma)=[\Gamma\sigma]\otimes \gamma\odot[\Gamma\sigma]^*_\gamma. \eqno(3)$$

Assume $\sigma_1,\sigma_2$ are elements of $G$ and that $[\Gamma\sigma_2\Gamma] = \bigcup [\Gamma\sigma_2 r_j]$.
Then from the Hecke algebra structure $[\Gamma\sigma_2\Gamma][\Gamma\sigma_1]=\sum_i[\Gamma\sigma_2 r_j\sigma_1]$.

Hence by using the identity (2)
\begin{align*}
\Phi(\Gamma\sigma_1\Gamma) \Phi(\Gamma\sigma_2) (\gamma) & =\chi_{\Gamma}([\Gamma\sigma_1\Gamma]\otimes[\Gamma\sigma_2\Gamma]^*)
([\Gamma\sigma_1]\otimes \gamma\odot [\Gamma\sigma_2]^*_\gamma)(e)
\\ & = \chi_{\Gamma}  \sum_{j,k} [\Gamma\sigma_1r_j \sigma_2]\otimes (\gamma \odot [\Gamma\sigma_1 r_k \sigma_2]^*_\gamma)(e).
\end{align*}

Passing $\chi_{\Gamma}$ to the right hand side, we obtain
$$
\sum_{j,k} [\Gamma\sigma_1r_j \sigma_2]\otimes (\gamma \odot [\Gamma\sigma_1 r_k \sigma_2]^*_\gamma)
\cdot \chi_{\sigma_1r_j \sigma_2 \Gamma \sigma_2^{-1} r_k^{-1} \sigma_1^{-1} \cap \Gamma}(e).
$$
Since in the expansion $[\Gamma\sigma_2\Gamma][\Gamma\sigma_1]=\sum_{j} [\Gamma\sigma_2r_j \sigma_1]$
to cosets $[\Gamma\sigma_2r_j \sigma_1]$ do not repeat themselves, the last sum is further equal to
$$
\sum_{j} [\Gamma\sigma_2r_j \sigma_1] \otimes 
(\gamma \odot [\Gamma\sigma_1 r_j \sigma_2]^*_\gamma)
$$
which is clearly equal by formula (3) to
$$
\Phi([\Gamma\sigma_2\Gamma] [\Gamma\sigma_1]) (\gamma)=\sum_{j}\Phi( [\Gamma\sigma_2r_j \sigma_1])(\gamma).
$$

\begin{rem}
The representation in the previous theorem is extended with the same method to 
${\rm Sp}\{[\sigma_1\Gamma\sigma_2]\mid \sigma_1,\sigma_2\in G\}$ by defining for
$\sigma_1,\sigma_2$ in $G$
$$
\Phi([\sigma_1\Gamma\sigma_2]) = \chi_\Gamma \bigg([\sigma_1\Gamma\sigma_2] \otimes
\sum_{[\Gamma z]\subseteq \Gamma\sigma_1\Gamma\sigma_2\Gamma} [\Gamma z]\bigg) \chi_\Gamma
$$
and
$$
\Phi(\sigma_1\Gamma) = \chi_\Gamma ([\sigma_1\Gamma]\otimes \Gamma\sigma_1\Gamma)\chi_\Gamma
$$
and verifying
$$
\Phi([\sigma_1\Gamma])\Phi([\Gamma\sigma_2]) = \Phi([\sigma_1\Gamma\sigma_2]).
$$
\end{rem}

We conclude by noting that in particular we obtained a canonical operator system structure on
${\rm Sp}\{[\sigma_1\Gamma][\Gamma\sigma_2]\mid \sigma_1,\sigma_2\in G\}$ 
with a trace
$$
\langle [\sigma_1\Gamma][\Gamma\sigma_2], [\sigma_3\Gamma][\Gamma\sigma_4]\rangle
= \tau \langle [\sigma_1\Gamma][\Gamma\sigma_2] ([\sigma_3\Gamma][\Gamma\sigma_4])^*\rangle.
$$

This scalar product has additional positivity properties that derive from ciclicity
of the trace ([Ra2]).

\begin{rem} 
Assume $G = \PSL_2(\Z[\frac{1}{p}])$, $\Gamma = \PSL_2(\Z)$. Then there exists a canonical isometry
$\Psi : \ell^2(\Gamma/G) \to \ell^2(F_N)$ where $N=\frac{p+1}{2}$, that is a 
$\C(\Gamma\backslash G/\Gamma) - A_{\rm red}$ bimodule map (where $A_{\rm red}$ is the algebra
generated in $\L(F_N)$ by $$\sum_{i=1}^N u_i + \sum_{i=1}^N u_i^{-1},$$ $u_i$ the generators of $F_N$ [Py]).

Then, since 
$${\rm Sp}[\sigma_1\Gamma][\Gamma\sigma_2] \cong \C(G/\Gamma) \textstyle\bigotimes\limits_{\C(\Gamma\backslash G/\Gamma)} \C(\Gamma/G)$$
it follows that the operator system structure corresponds to an operator system structure on
$\ell^2(F_N)\bigotimes\limits_A \ell^2 (F_N)$.
\end{rem}

{\it Proof.} Indeed to define $\Psi$ one devides the cosets $[\Gamma\sigma_p s_i]_{i=1}^{p+1}$,
$\sigma_p = \left(\begin{array}{cc} p & 0 \\ 0 & 1 \end{array} \right)$ into two
parts, one which is mapped into $u_1,\ldots,u_N$ the other into $u_1^{-1},\ldots,u_N^{-1}$.~$\quad\Box$

\begin{rem} 
Note that there exists a canonical embedding $$\ell^2 (G/\Gamma_\sigma)\textstyle\bigotimes\limits_{\C(\Gamma_\sigma\backslash G/\Gamma_\sigma)} \ell^2 (\Gamma_\sigma\backslash G)$$
into 
$\ell^2 (G/\Gamma_{\sigma_1})\bigotimes\limits_{\C(\Gamma_{\sigma_1}\backslash G/\Gamma_{\sigma_1})} \ell^2 (\Gamma_{\sigma_1}\backslash G)$ if $\Gamma_{\sigma_1} \subseteq \Gamma_{\sigma_0}$.

Indeed if the identity (1) is satisfied for $\sigma_1^i, \sigma_2^i, \theta_1^j,\theta_2^j$ with $\Gamma_{\sigma_1}$
instead of $\Gamma$, it will be satisfied also for $\Gamma_{\sigma_0}$ instead of $\Gamma_{\sigma_1}$
(this is proven by splitting $\Gamma_{\sigma_0}$ into $\Gamma_{\sigma_1}$ cosets). Hence one can define (in
analogy with the Hecke modular algebra from ([CoMo]) an operator system structure on the inductive limit
$$\textstyle\bigcup\limits_{\Gamma_\sigma \in\S}
\C(G\backslash\Gamma_\sigma)\bigotimes\limits_{\C(\Gamma_\sigma\backslash G/\Gamma_\sigma)} 
\C(\Gamma_\sigma\backslash G).$$
\end{rem}

\end{document}